\documentclass{amsproc}%
\usepackage{amssymb}
\usepackage{graphicx}
\usepackage{amscd}
\usepackage{amsmath}
\usepackage{amsfonts}%
\theoremstyle{plain}

\numberwithin{equation}{section}

\begin{document}
\title[Which self-maps appear as lattice endomorphisms? ]{Which self-maps appear as lattice endomorphisms?}
\author{Jen\H{o} Szigeti}
\address{Institute of Mathematics, University of Miskolc, Miskolc, Hungary 3515}
\email{matjeno@uni-miskolc.hu}
\thanks{The author was supported by OTKA K101515 and his research was carried out as
part of the TAMOP-4.2.1.B-10/2/KONV-2010-0001 project with support by the
European Union, co-financed by the European Social Fund.}

\begin{abstract}
Let $f:A\rightarrow A$ be a self-map of the set $A$. We give a necessary and
sufficient condition for the existence of a lattice structure $(A,\vee
,\wedge)$ on $A$ such that $f$ becomes a lattice endomorphism with respect to
this structure.

\end{abstract}
\maketitle

\noindent1. INTRODUCTION

\bigskip

A \textit{partially ordered set (poset)} is a set $P$\ together with a
reflexive, antisymmetric, and transitive (binary) relation $r\subseteq P\times
P$. For $(x,y)\in r$ we write $x\leq_{r}y$ or simply $x\leq y$. If $r\subseteq
r^{\prime}\subseteq P\times P$ for the partial orders $r$ and $r^{\prime}$,
then $r^{\prime}$ is an \textit{extension} of $r$. A map $p:P\longrightarrow
P$ is \textit{order-preserving} if $x\leq y$ implies $p(x)\leq p(y)$ for all
$x,y\in P$. The poset $(P,\leq)$ is a \textit{lattice} if any two elements
$x,y\in P$ have a unique least upper bound (lub) $x\vee y$ and a unique
greatest lower bound (glb) $x\wedge y$ (in $P$). The operations $\vee$ and
$\wedge$ are associative, commutative, and satisfy the following absorption
laws: $(x\vee y)\wedge y=y$ and $(x\wedge y)\vee y=y$. Any binary operations
$\vee$ and $\wedge$ on $P$ having these properties define a binary relation
$r=\{(x,x\vee y):x,y\in P\}\subseteq P\times P$ on $P$, which is a partial
order. In fact $(P,\leq_{r})$ is a lattice with lub $\vee$ and glb $\wedge$.
Lattices play a fundamental role in many areas of mathematics (see [1],[3]).

In the present paper we consider a self-map $f:A\longrightarrow A$\ of a set
$A$. A list $x_{1},\ldots,x_{n}$ of distinct elements from $A$ is a
\textit{cycle} (of \textit{length} $n$) with respect to $f$ if $f(x_{i}%
)=x_{i+1}$ for each $1\leq i\leq n-1$ and also $f(x_{n})=x_{1}$. A
\textit{fixed point} of the function $f$ is a cycle of length $1$, i.e. an
element $x_{1}\in A$ with $f(x_{1})=x_{1}$. A cycle that is not a fixed point
is \textit{proper}.

If $(A,\vee,\wedge)$ is a lattice (on the set $A$) such that $f(x\vee
y)=f(x)\vee f(y)$ and $f(x\wedge y)=f(x)\wedge f(y)$ for all $x,y\in A$, then
$f$ is a \textit{lattice endomorphism} of $(A,\vee,\wedge)$. A lattice
endomorphism is an order-preserving map (with respect to the order relation of
the lattice), but the converse is not true in general. For a proper cycle
$x_{1},\ldots,x_{n}\in A$ with respect to a lattice endomorphism $f$, if we
put%
\[
p=x_{1}\vee x_{2}\vee\cdots\vee x_{n}%
\]
and%
\[
q=x_{1}\wedge x_{2}\wedge\cdots\wedge x_{n},
\]
then $p\neq q$. The equalities%
\[
f(p)=f(x_{1})\vee f(x_{2})\vee\cdots\vee f(x_{n})=x_{2}\vee\cdots\vee
x_{n}\vee x_{1}=p,
\]%
\[
f(q)=f(x_{1})\wedge f(x_{2})\wedge\cdots\wedge f(x_{n})=x_{2}\wedge
\cdots\wedge x_{n}\wedge x_{1}=q
\]
show that $p$ and $q$ are distinct fixed points of $f$. It follows that any
lattice endomorphism having a proper cycle must have at least two fixed points.

We prove that the above combinatorial property completely characterizes the
possible lattice endomorphisms. More precisely, for a map $f:A\longrightarrow
A$\ there exists a lattice $(A,\vee,\wedge)$ on $A$ such that $f$ is a lattice
endomorphism of $(A,\vee,\wedge)$ if and only if $f$ has no proper cycles or
$f$ has at least two fixed points.

The construction in the proof of our main result is based on the use of the
maximal $f$-compatible extensions of an $f$-compatible partial order relation
on $A$. Such extensions were completely determined in [2] and [5]. In order to
make the exposition more self-contained, we present the necessary background
about maximal compatible extensions.

\bigskip

\noindent2. PRELIMINARY DEFINITIONS AND\ RESULTS

\bigskip

\noindent Let $f:A\longrightarrow A$ be a function, and define the equivalence
relation $\sim_{f}$ as follows: for $x,y\in A$, let $x\sim_{f}y$ if
$f^{k}(x)=f^{l}(y)$ for some integers $k\geq0$ and $l\geq0$. The equivalence
class $\left[  x\right]  _{f}$\ of an element $x\in A$ is the $f$%
\textit{-component} of $x$. We note that $\left[  x\right]  _{f}$ is closed
with respect to the action of $f$ and hence contains the $f$\textit{-orbit}
of$~x$:%
\[
\{x,f(x),\ldots,f^{k}(x),\ldots\}\subseteq\left[  x\right]  _{f}\text{ }.
\]
An element $c\in A$ is \textit{cyclic} with respect to $f$ if $f^{m}(c)=c$ for
some integer $m\geq1$. The \textit{period} of a cyclic element $c$, written as
$n(c)$, is defined by%
\[
n(c)=\min\{m:m\geq1\text{ and }f^{m}(c)=c\},
\]
and $f^{k}(c)=f^{l}(c)$ holds if and only if $k-l$ is divisible by $n$. The
\textit{full cycle }of a cyclic element $c$ is the $f$-orbit
$\{c,f(c),...,f^{n(c)-1}(c)\}$. The $f$-orbit of $x$ is finite if and only if
$\left[  x\right]  _{f}$\ contains a cyclic element. If $c_{1},c_{2}\in\left[
x\right]  _{f}$ are cyclic elements, then $n(c_{1})=n(c_{2})=n(x)$, and this
number is the \textit{period }of $x$. If the $f$-orbit of $x$ is infinite,
then put $n(x)=\infty$. Clearly, $x\sim_{f}y$ implies $n(x)=n(y)$. We note
that the presence of a cyclic element in $\left[  x\right]  _{f}$ does not
imply that $\left[  x\right]  _{f}$ is finite. The function $f$ has a proper
cycle if there exists a cyclic element $c\in A$ with respect to $f$ such that
$n(c)\geq2$.

\bigskip

\noindent\textbf{2.1. Theorem (see [4]).}\textit{\ Let }$r$\textit{\ be an
order relation on the set }$A$\textit{, and let }$f:A\rightarrow
A$\textit{\ be an order-preserving map with respect to }$r$\textit{. If there
is no proper cycle of }$f$\textit{, then there exists a linear extension }%
$R$\textit{ of }$r$\textit{ such that }$f$\textit{ is order-preserving with
respect to }$R$.

\bigskip

\noindent\textbf{2.2. Corollary.} \textit{If }$f:A\longrightarrow
A$\textit{\ is a function with no proper cycles, then there exists a
distributive lattice }$(A,\vee,\wedge)$\textit{ on }$A$\textit{ such that }%
$f$\textit{ is a lattice endomorphism of }$(A,\vee,\wedge)$.

\bigskip

\noindent The following definitions appear in [2]. A pair $(x,y)\in A\times A$
is $f$\textit{-prohibited} if there exist integers $k$, $l$, and $n$ with
$k\geq0$, $l\geq0$, and $n\geq2$ such that $n$ is not a divisor of $k-l$, the
elements $f^{k}(x),f^{k+1}(x),...,f^{k+n-1}(x)$ are distinct and
$f^{k+n}(x)=f^{k}(x)=f^{l}(y)$. For an $f$-prohibited pair $(x,y)$ and
integers $k$ and $n$\ as above, $f^{k}(x)$ is a cyclic element in $\left[
x\right]  _{f}=\left[  y\right]  _{f}$\ of period $n$. The \textit{distance}
$d(y,c)$ between an element $y\in\left[  x\right]  _{f}$ and a given cyclic
element $c\in\left[  x\right]  _{f}$\ (of period $n\geq1$) is defined by%
\[
d(y,c)=\min\{t:t\geq0\text{\ and }f^{t}(y)=c\}.
\]
Clearly, $f^{t}(y)=c$\ holds if and only if $t\geq d(y,c)$\ and $t-d(y,c)$\ is
divisible by $n$. We note that $d(f(c),c)=n-1$,\ and for $y\neq c$\ we have
$d(f(y),c)=d(y,c)-1$. It is straightforward to see that $(x,y)$\ is
$f$-prohibited if and only if $d(x,c)-d(y,c)$\ is not divisible by $n$.

\bigskip

\noindent\textbf{2.3. Proposition (see [2]).}\textit{\ Let }$r$\textit{\ be an
order relation on the set }$A$\textit{\ and }$f:A\rightarrow A$\textit{\ be an
order-preserving map with respect to }$r$. \textit{If }$(x,y)\in A\times
A$\textit{\ is an }$f$\textit{-prohibited pair, then }$x$\textit{\ and }%
$y$\textit{\ are incomparable with respect to }$r$\textit{.}

\bigskip

\noindent\textbf{2.4. Lemma (see [2]).}\textit{\ Let }$f:A\longrightarrow
A$\textit{ be a self-map on a set }$A$\textit{. Let }$c\mathit{\ }$\textit{be
a cyclic element, with }$c\in\left[  x\right]  _{f}$\textit{\ for some }$x\in
A$\textit{. If }$r$\textit{\ is an order relation on }$A$\textit{, and }%
$f$\textit{ is order-preserving with respect to }$r$\textit{, then there
exists an order relation }$\rho$\textit{\ on }$\left[  x\right]  _{f}%
$\textit{\ with the following properties:}

\noindent1. $\rho$\textit{\ is an extension of }$r$\textit{ (on }$\left[
x\right]  _{f}$\textit{): }$r\cap(\left[  x\right]  _{f}\times\left[
x\right]  _{f})\subseteq\rho$\textit{,}

\noindent2. $f$\textit{ is order-preserving with respect to }$\rho$\textit{,}

\noindent3. $\left[  x\right]  _{f}$\textit{\ is the disjoint union of sets}
$E_{0},\ldots,E_{n-1}$\textit{ and each}%
\[
E_{i}=\{u\in\left[  x\right]  _{f}:d(u,c)-i\text{ is divisible by
}n(c)\}\text{ },\text{ }0\leq i\leq n-1
\]

\noindent\ \ \ \ \textit{is a chain with respect to }$\rho$\textit{,}

\noindent4.\textit{ }$f(E_{0})\subseteq E_{n-1}$\textit{ and }$f(E_{i}%
)\subseteq E_{i-1}$\textit{\ for\ }$1\leq i\leq n-1$\textit{,}

\noindent5. \textit{any element }$(u,v)\in E_{i}\times E_{j}$\textit{\ with
}$i\neq j$\textit{ is an }$f$\textit{-prohibited pair, and the set}

\noindent\ \ \ \ $\{u,v\}$ \textit{has no upper and lower bounds in }$\left[
x\right]  _{f}$\textit{\ with respect to }$\rho$\textit{.}

\bigskip

\noindent3. MAKING\thinspace$f$\thinspace A\thinspace LATTICE\thinspace ENDOMORPHISM

\bigskip

\noindent\textbf{3.1. Theorem.}\textit{\ Let }$r$\textit{\ be an order
relation on the set }$A$\textit{, and let }$f:A\longrightarrow A$\textit{\ be
an order-preserving map with respect to }$r$\textit{ having distinct fixed
points }$p,q\in A$\textit{. If }$x$\textit{ and }$y$\textit{ are }%
$r$\textit{-incomparable for all }$x,y\in A$\textit{ such that }$\left[
x\right]  _{f}\neq\left[  y\right]  _{f}$\textit{ and }$2\leq n(x)\neq\infty
$\textit{, then there exists an extension }$R$\textit{ of }$r$\textit{ such
that }$(A,\leq_{R})$\textit{ is a lattice and }$f$\textit{ is a lattice
endomorphism of }$(A,\leq_{R})$\textit{.}

\bigskip

\noindent\textbf{Proof.} Let%
\[
A_{0}=\{x\in A:\left[  x\right]  _{f}\text{ contains a proper cycle}\}=\{x\in
A:2\leq n(x)\neq\infty\}
\]
The set $A_{0}$\ is the $f$-cyclic part of $A$. Let%
\[
A_{\ast}\!=\!A\setminus A_{0}\!=\!\{x\in A:\left[  x\right]  _{f}\text{ has no
proper cycle}\}\!=\!\{x\in A:n(x)=1\text{ or }n(x)=\infty\}
\]
The set $A_{\ast}$ is the $f$-acyclic part of $A$. We have either $\left[
x\right]  _{f}\subseteq A_{0}$ or $\left[  x\right]  _{f}\subseteq A_{\ast}$
for all $x\in A$. Clearly, both $A_{0}$ and $A_{\ast}$ are closed with respect
to the action of $f$, i.e. $f(A_{0})\subseteq A_{0}$ and $f(A_{\ast})\subseteq
A_{\ast}$. Since $f:A_{\ast}\longrightarrow A_{\ast}$ has no proper cycle (in
$A_{\ast}$), Theorem 2.1 ensures the existence of a linear extension $R_{\ast
}$\ of $r\cap(A_{\ast}\times A_{\ast})$ (on $A_{\ast}$) such that $f$ is
order-preserving with respect to $R_{\ast}$. In view of $p,q\in A_{\ast}$, we
may assume $p\leq_{R_{\ast}}q$.

For an appropriate subset $\{x_{t}:t\in T\}$ of $A_{0}$, where the indices are
taken from an idex set $T$, we have $\{\left[  x\right]  _{f}:x\in
A_{0}\}=\{\left[  x_{t}\right]  _{f}:t\in T\}$, and $\left[  x_{t}\right]
_{f}\neq\left[  x_{s}\right]  _{f}$ for all $t,s\in T$ with $t\neq s$. Such a
subset $\{x_{t}:t\in T\}\subseteq A_{0}$ is an \textit{irredundant set of
representatives} of the equivalence classes of $\sim_{f}$ (in $A_{0}$). That
is%
\[
A_{0}=%
%TCIMACRO{\tbigcup \nolimits_{t\in T}}%
%BeginExpansion
{\textstyle\bigcup\nolimits_{t\in T}}
%EndExpansion
\left[  x_{t}\right]  _{f}\text{ and }\left[  x_{t}\right]  _{f}\cap\left[
x_{s}\right]  _{f}=\varnothing\text{ for all }t,s\in T\text{ with }t\neq s.
\]
Call two elements of $A$ \textit{concurrent} if some power of $f$ maps them to
the same element. Concurrency is an equivalence relation finer than $\sim_{f}%
$. For $t\in T$, the $\sim_{f}$-class of $x_{t}$ is partitioned into
$n(x_{t})$ concurrency classes:
\[
\left[  x_{t}\right]  _{f}=E_{0}^{(t)}\cup E_{1}^{(t)}\cup...\cup
E_{n(x_{t})-1}^{(t)}\text{, where}%
\]%
\[
E_{i}^{(t)}=\{u\in\left[  x_{t}\right]  _{f}:d(u,c)-i\text{ is divisible by
}n(x_{t})\}
\]
for some fixed cyclic element $c\in\left[  x_{t}\right]  _{f}$. Application of
Lemma 2.4 gives the existence of a partial order extension $\rho_{t}$ of $r$
on $\left[  x_{t}\right]  _{f}$ ($r\cap(\left[  x_{t}\right]  _{f}%
\times\left[  x_{t}\right]  _{f})\subseteq\rho_{t}$ holds) such that $f$
preserves $\rho_{t}$ and each $E_{i}^{(t)}$ is a chain with respect to
$\rho_{t}$.

Take the following subsets of $A\times A$:%
\[
P\!=\!\{\!(a,x)\!:\!a\!\in\!A_{\ast},x\!\in\!A_{0}\text{ and }a\!\leq
_{R_{\ast}}\!p\}\text{ and }Q\!=\!\{\!(y,b)\!:\!b\!\in\!A_{\ast},y\!\in
\!A_{0}\text{ and }q\!\leq_{R_{\ast}}\!b\}.
\]

Let%
\[
R=R_{\ast}\cup\left(
%TCIMACRO{\tbigcup \nolimits_{t\in T}}%
%BeginExpansion
{\textstyle\bigcup\nolimits_{t\in T}}
%EndExpansion
\rho_{t}\right)  \cup P\cup Q.
\]
We claim that $R$ is an extension of $r$\ that is a lattice and that $f$ is a
lattice endomorphism of $(A,\leq_{R},\vee,\wedge)$. The proof consists of the
following straightforward steps.

Notice that $R_{\ast}\subseteq A_{\ast}\times A_{\ast}$, $\rho_{t}%
\subseteq\left[  x_{t}\right]  _{f}\times\left[  x_{t}\right]  _{f}\subseteq
A_{0}\times A_{0}$, $P\subseteq A_{\ast}\times A_{0}$, and $Q\subseteq
A_{0}\times A_{\ast}$. Also the direct products $A_{\ast}\times A_{\ast}$,
$A_{\ast}\times A_{0}$, $A_{0}\times A_{\ast}$, and $\left[  x_{t}\right]
_{f}\times\left[  x_{t}\right]  _{f}$ (for $t\in T$) are pairwise disjoint.

In order to see $r\subseteq R$, take $(u,v)\in r$.

\noindent(1) If $(u,v)\in A_{\ast}\times A_{\ast}$, then $r\cap(A_{\ast}\times
A_{\ast})\subseteq R_{\ast}$ implies $(u,v)\in R$.

\noindent(2) If $(u,v)\in A_{\ast}\times A_{0}$, then $\left[  u\right]
_{f}\neq\left[  v\right]  _{f}$ and $2\leq n(v)\neq\infty$ contradicts
$(u,v)\in r$.

\noindent(3) $(u,v)\in A_{0}\times A_{\ast}$ is also impossible.

\noindent(4) If $(u,v)\in A_{0}\times A_{0}$, then $(u,v)\in\left[
x_{t}\right]  _{f}\times\left[  x_{s}\right]  _{f}$ for some $t,s\in T$.
Clearly, $t\neq s$ would imply $\left[  u\right]  _{f}\neq\left[  v\right]
_{f}$, and then $2\leq n(u)\neq\infty$ contradicts $(u,v)\in r$. Thus $t=s$,
and $r\cap(\left[  x_{t}\right]  _{f}\times\left[  x_{t}\right]
_{f})\subseteq\rho_{t}$ yields $(u,v)\in R$.

We prove that $R$ is a partial order.

\noindent Antisymmetry: Let $(u,v)\in R$ and $(v,u)\in R$.

\noindent(1) If $(u,v),(v,u)\in R_{\ast}$, then $u=v$ follows from the
antisymmetric property of $R_{\ast}$

\noindent(2) If $(u,v)\in\rho_{t}$ and $(v,u)\in\rho_{s}$, then $t=s$, and
$u=v$ follows from the antisymmetric property of $\rho_{t}$.

\noindent(3) If $(u,v)\in P$ and $(v,u)\in Q$, then $u\leq_{R_{\ast}}p$ and
$q\leq_{R_{\ast}}u$ imply $q\leq_{R_{\ast}}p$, contradicting with
$p\leq_{R_{\ast}}q$ and $p\neq q$.

\noindent(4) If $(u,v)\in Q$ and $(v,u)\in P$, then interchanging the roles of
$u$ and $v$ leads to a similar contradiction as in case (3).

\noindent Transitivity: Let $(u,v)\in R$ and $(v,w)\in R$.

\noindent(1) If $(u,v),(v,w)\in R_{\ast}$, then $(u,w)\in R_{\ast}$ follows
from the transitivity of $R_{\ast}$.

\noindent(2) If $(u,v)\in R_{\ast}$ and $(v,w)\in P$, then $u\leq_{R_{\ast}%
}v\leq_{R_{\ast}}p$ and $w\in A_{0}$ imply $(u,w)\in P$.

\noindent(3) If $(u,v)\in\rho_{t}$ and $(v,w)\in\rho_{s}$, then we have $t=s$,
and $(u,w)\in\rho_{t}$ follows from the transitivity of $\rho_{t}$.

\noindent(4) If $(u,v)\in\rho_{t}$ and $(v,w)\in Q$, then $u,v\in A_{0}$,
$w\in A_{\ast}$, and $q\leq_{R_{\ast}}w$. It follows that $(u,w)\in Q$.

\noindent(5) If $(u,v)\in P$ and $(v,w)\in\rho_{t}$, then $v,w\in A_{0}$,
$u\in A_{\ast}$, and $u\leq_{R_{\ast}}p$. It follows that $(u,w)\in P$.

\noindent(6) If $(u,v)\in P$ and $(v,w)\in Q$, then $u\leq_{R_{\ast}}%
p\leq_{R_{\ast}}q\leq_{R_{\ast}}w$, from which $(u,w)\in R_{\ast}$ follows.

\noindent(7) If $(u,v)\in Q$ and $(v,w)\in R_{\ast}$, then $u\in A_{0}$ and
$q\leq_{R_{\ast}}v\leq_{R_{\ast}}w$ imply $(u,w)\in P$.

\noindent(8) If $(u,v)\in Q$ and $(v,w)\in P$, then $q\leq_{R_{\ast}}%
v\leq_{R_{\ast}}p$ contradicts $p\leq_{R_{\ast}}q$ and $p\neq q$.

We note that $f$\ is order-preserving with respect to $(A_{\ast},\leq
_{R_{\ast}})$, and $(\left[  x_{t}\right]  _{f},\rho_{t})$ for $t\in T$. In
order to check the order-preserving property of $f$ with respect to
$(A,\leq_{R})$, it is enough to see that $(a,x)\in P$ implies $(f(a),f(x))\in
P$ and $(y,b)\in Q$ implies $(f(y),f(b))\in Q$. Obviously, $a\in A_{\ast}$,
$x\in A_{0}$, and $a\leq_{R_{\ast}}p$ imply $f(a)\in A_{\ast}$, $f(x)\in
A_{0}$, and $f(a)\leq_{R_{\ast}}f(p)=p$. Similarly, $b\in A_{\ast}$, $y\in
A_{0}$, and $q\leq_{R_{\ast}}b$ imply $f(b)\in A_{\ast}$, $f(y)\in A_{0}$, and
$q=f(q)\leq_{R_{\ast}}f(b)$.

If $u,v\in A$ are comparable elements with respect to $R$, then the existence
of the supremum $u\vee v$ and the infimum $u\wedge v$ in $(A,\leq_{R})$ is
evident; moreover, the order-preserving property of $f$ ensures that%
\[
f(u\vee v)=f(u)\vee f(v)\text{ and }f(u\wedge v)=f(u)\wedge f(v).
\]

If $u,v\in A$ are incomparable elements with respect to $R$, then we have the
following possibilities.

\noindent(1) If $u\in A_{\ast}$ and $v\in A_{0}$, then $(u,v)\notin P$,
$(v,u)\notin Q$, and the linearity of $R_{\ast}$\ imply $p\leq_{R_{\ast}}%
u\leq_{R_{\ast}}q$, from which $u\vee v=q$ and $u\wedge v=p$ follow in
$(A,\leq_{R})$. Since $f(u)\in A_{\ast}$, $f(v)\in A_{0}$, and $p=f(p)\leq
_{R_{\ast}}f(u)\leq_{R_{\ast}}f(q)=q$, we deduce that%
\[
f(u\vee v)=f(q)=q=f(u)\vee f(v)\text{ and }f(u\wedge v)=f(p)=p=f(u)\wedge
f(v).
\]
\noindent(2) If $u\in A_{0}$ and $v\in A_{\ast}$, then interchanging the roles
of $u$ and $v$ leads to the same result as in case (1).

\noindent(3) If $u,v\in A_{0}$ and $\left[  u\right]  _{f}\neq\left[
v\right]  _{f}$, then $u\vee v=q$ and $u\wedge v=p$ in $(A,\leq_{R})$ follow
directly from the definition of $R$. Since $f(u),f(v)\in A_{0}$, and $\left[
f(u)\right]  _{f}=\left[  u\right]  _{f}\neq\left[  v\right]  _{f}=\left[
f(v)\right]  _{f}$, we deduce%
\[
f(u\vee v)=f(q)=q=f(u)\vee f(v)\text{ and }f(u\wedge v)=f(p)=p=f(u)\wedge
f(v).
\]
\noindent(4) If $u,v\in A_{0}$ and $\left[  u\right]  _{f}=\left[  v\right]
_{f}=\left[  x_{t}\right]  _{f}$ for some unique $t\in T$, then $(u,v)\notin
\rho_{t}$ implies $(u,v)\in E_{i}^{(t)}\times E_{j}^{(t)}$\ for some unique
$0\leq i,j\leq n(x_{t})-1$ with $i\neq j$. In view of $E_{i}^{(t)}\cap
E_{j}^{(t)}=\varnothing$ and (5) of Lemma 2.4, we conclude that the set
$\{u,v\}$ has no upper and lower bounds in $(\left[  x_{t}\right]  _{f}%
,\rho_{t})$. It follows that $u\vee v=q$ and $u\wedge v=p$ in $(A,\leq_{R})$.

Since $f(E_{i})\subseteq E_{i-1}$ implies $(f(u),f(v))\in E_{i-1}^{(t)}\times
E_{j-1}^{(t)}$\ (notice that $E_{-1}^{(t)}=E_{n(x_{t})-1}^{(t)}$), we deduce
in a similar way%
\[
f(u)\vee f(v)=q=f(q)=f(u\vee v)\text{ and }f(u)\wedge f(v)=p=f(p)=f(u\wedge
v).\text{ }\square
\]

\bigskip

\noindent\textbf{3.2. Corollary.} \textit{If the number of fixed points of the
function }$f:A\longrightarrow A$\textit{\ is at least }$2$\textit{, then there
exists a lattice structure }$(A,\vee,\wedge)$\textit{ on }$A$\textit{ such
that }$f$\textit{ is a lattice endomorphism of }$(A,\vee,\wedge)$.

\bigskip

\noindent\textbf{Proof.} Let $p$ and $q$ be distinct fixed points of $f$. The
application of Theorem 3.1 yields a partial order extension $R$ of the
identity partial order $\{(x,x):x\in A\}$ such that $(A,\leq_{R},\vee,\wedge)$
is a lattice and $f$ is a lattice endomorphism of $(A,\leq_{R},\vee,\wedge)$.
$\square$

\bigskip

\noindent The combination of Corollaries 2.2 and 3.2 provides the complete
answer (formulated in the introduction) to the question in the title of the
paper. We pose a further problem.

\bigskip

\noindent\textbf{3.3. Problem.}\textit{ Consider an arbitrary function
}$f:A\longrightarrow A$\textit{. Find necessary and sufficient conditions for
the existence of a modular (or distributive) lattice structure }%
$(A,\vee,\wedge)$\textit{ on }$A$\textit{ such that }$f$\textit{ becomes a
lattice endomorphism of }$(A,\vee,\wedge)$\textit{. The similar question seems
to be interesting for other algebraic structures such as (Abelian) groups,
rings and modules.}

\bigskip

\noindent\textbf{3.4. Example.} Let $A=\{p,q,x_{1},x_{2},\ldots,x_{n}\}$,
where $n\geq3$, and let $f:A\longrightarrow A$ be a function with $f(p)=p$,
$f(q)=q$, $f(x_{n})=x_{1}$, and $f(x_{i})=x_{i+1}$ for $1\leq i\leq n-1$. If
$f$ is an endomorphism of some lattice $(A,\leq,\vee,\wedge)$, then $f$ is
order-preserving with respect to $(A,\leq)$, and Proposition 2.3 ensures that
the proper cycle $\{x_{1},\ldots,x_{n}\}$ of $f$\ is an antichain in
$(A,\leq)$. Since $x_{1}\vee\cdots\vee x_{n}$ and $x_{1}\wedge\cdots\wedge
x_{n}$ are distinct fixed points of $f$, one of $x_{1}\vee\cdots\vee x_{n}$
and $x_{1}\wedge\cdots\wedge x_{n}$ is $p$ and the other is $q$. Thus
$(A,\leq,\vee,\wedge)$ is isomorphic to the lattice $M_{n}$ in both cases. It
follows that there is no distributive lattice structure on $A$ making $f$ a
lattice endomorphism (even though $f$ has two fixed points).

\bigskip

\noindent REFERENCES

\bigskip

\noindent\lbrack1] Birkhoff, G.\textit{ Lattice Theory,} Colloquium
Publications, Vol. 25, AMS Providence (Third Edition, 1995)

\noindent\lbrack2] Foldes, S.; Szigeti, J. \textit{Maximal compatible
extensions of partial orders, }J. Australian Math. Soc. 81 (2006), 245-252.

\noindent\lbrack3] Gratzer, G. \textit{General Lattice Theory,} Birkhauser
Verlag, Basel-Boston-Berlin (2003)

\noindent\lbrack4] Szigeti, J. , Nagy, B. \textit{Linear extensions of partial
orders preserving monotonicity,} Order 4 (1987), 31-35.

\noindent\lbrack5] Szil\'{a}gyi, Sz. \textit{A characterization and the
intersection of the maximal compatible extensions of a partial order,} Order
Vol. 25, 4 (2008), 321-334.

\end{document}